\documentclass[runningheads]{llncs}
\usepackage{graphicx} 
\usepackage[T1]{fontenc}
\usepackage{newunicodechar}
\newunicodechar{ℓ}{\ensuremath{\ell}}
\usepackage{listings}
\usepackage{amsmath}
\usepackage{amssymb}

\usepackage{color}
\definecolor{keywordcolor}{rgb}{0.7, 0.1, 0.1}   
\definecolor{tacticcolor}{rgb}{0.0, 0.1, 0.6}    
\definecolor{commentcolor}{rgb}{0.4, 0.4, 0.4}   
\definecolor{symbolcolor}{rgb}{0.0, 0.1, 0.6}    
\definecolor{sortcolor}{rgb}{0.1, 0.5, 0.1}      
\definecolor{attributecolor}{rgb}{0.7, 0.1, 0.1}

\title{Formalizing $A_1^{(1)}$ Curve Neighborhoods in Lean 4}
\titlerunning{Formalizing $A_1^{(1)}$ Curve Neighborhoods in Lean 4}
\author{Yihe Huang\inst{1,2} \and Sizhe Cui\inst{2} \and Jiaqi Wang\inst{2} \and Jujian Zhang\inst{3}}
\institute{
Department of Mathematics, University of Chinese Academy of Sciences, Beijing 100049, China\\
\email{huangyihe22@mails.ucas.ac.cn}
\and
State Key Laboratory of Mathematical Sciences, Academy of Mathematics and Systems Science, University of Chinese Academy of Sciences, Beijing 100190, China\\
\email{\{wangjiaqi21c,cuisizhe25\}@mails.ucas.ac.cn}
\and
Department of Mathematics, Imperial College London, London SW7 2AZ, United Kingdom\\
\email{jujian.zhang19@imperial.ac.uk}
}
\date{January 2026}

\begin{document}

\maketitle

\begin{abstract}
Combinatorial curve neighborhoods are somewhat foundational when setting up the quantum Schubert calculus for affine flag manifolds. In the specific case of type $A_1^{(1)}$, you can encode these neighborhoods entirely within the moment graph of the infinite dihedral group $D_\infty$. Building on the framework developed by Mihalcea and Norton, this paper presents a complete, axiom-free formalization of these combinatorial curve neighborhoods in Lean 4. Rather than just wrapping mathematical statements, we formalized $D_\infty$ directly as a Coxeter system to explicitly compute length functions and degree maps. Reachable sets are defined through edge chains bounded by specific degrees, and we ultimately characterize the curve neighborhood by the maximal vertices inside these sets. The core effort here lies in formally verifying the explicit combinatorial formulas for curve neighborhoods of arbitrary elements. Interestingly, by restricting our search space to finite sets, we also managed to extract a fully computable version of these neighborhoods. 
\end{abstract}

\section{Introduction}
When analyzing the quantum cohomology and quantum $K$-theory rings of generalized flag manifolds, one often leans heavily on geometric \emph{curve neighborhoods} \cite{buch2015}. Mare and Mihalcea later extended this idea to affine variants \cite{mare2014}. To actually calculate these neighborhoods, it helps to encode them into the \emph{moment graph} of the affine flag manifold. Now, working with arbitrary Kac-Moody groups can get exceedingly complicated quickly. However, for type $A_1^{(1)}$—where the affine Weyl group is just the infinite dihedral group $D_\infty$—Mihalcea and Norton managed to provide a very explicit formula \cite{mihalcea2017}

Their proofs, while elegant, involve fairly intricate combinatorial steps: checking parity, verifying word lengths over $D_\infty$, and so on. These kinds of arguments are notoriously easy to slip up on by hand, making them a prime target for formal verification. In this paper, we present a complete, axiom-free formalization of combinatorial curve neighborhoods for type $A_1^{(1)}$. Using Lean 4~\cite{lean4} and Mathlib~\cite{mathlib}, we encode $D_\infty$ as a Coxeter system and formalize the underlying moment graph, ultimately securing a machine-checked proof of the main formula driving these curve neighborhoods.

\section{Mathematical Background}
\label{sec:background}

We briefly recall definitions following Mihalcea and Norton \cite{mihalcea2017}, forming the blueprint for our Lean 4 formalization.

\subsection{Dihedral Group and Moment Graph}
The infinite dihedral group $D_\infty$, the affine Weyl group of type $A_1^{(1)}$, is the Coxeter group generated by $s_0$ and $s_1$ subject to:
\begin{equation}
    s_0^2 = s_1^2 = 1.
\end{equation}
Elements $w \in D_\infty$ decompose uniquely into alternating products, whose \emph{length} $\ell(w)$ counts the reflections. Elements are classified into even-length \textbf{Rotations} $r(k)$ and odd-length \textbf{Reflections} $sr(k)$, mirroring constructors of \lstinline{DihedralGroup 0} in Mathlib. A root $\alpha = (a, b) \in \mathbb{Z}_{\ge 0}^2$ requires $|a - b| = 1$, corresponding to a root reflection $s_\alpha$ of length $a + b$.

The \emph{moment graph} $G = (V, E)$ has a vertex set $V = D_\infty$. A directed edge $u \xrightarrow{\alpha} v$ of \emph{degree} $\alpha$ exists if and only if $v = u s_\alpha$ for some root reflection $s_\alpha$:
\begin{equation}
    v = u s_\alpha.
\end{equation}
Let $\varphi: D_\infty \to \mathbb{Z}_{\ge 0}^2$ map $w$ to counts $(d_0, d_1)$ of $s_0$ and $s_1$ in its reduced word. Degrees are partially ordered component-wise.
An \emph{increasing chain} from $u_0$ to $u_n$ strictly increases in length at each step. Its degree is the sum of its edge degrees. Here, the Bruhat order $u < v$ coincides with an increasing chain, and $u < v \iff \ell(u) < \ell(v)$ \cite[Lemma 2.3]{mihalcea2017}.

\subsection{Combinatorial Curve Neighborhoods}
We are now ready to define the central mathematical object of our study.

\begin{definition}[Combinatorial Curve Neighborhood]
Fix a degree $d = (d_0, d_1)$ and an element $u \in D_\infty$. The combinatorial curve neighborhood, denoted $\Gamma_d(u)$, is the set of elements $v \in D_\infty$ satisfying:
\begin{enumerate}
    \item There exists a chain of some degree $d' \le d$ from $u$ to $v$ in the moment graph $G$.
    \item The elements $v$ are maximal (with respect to the Bruhat order) among all elements satisfying condition (1).
\end{enumerate}
\end{definition}

Computing $\Gamma_d(u)$ directly from its graph definition turns out to be highly non-trivial, simply because the chains are recursive by nature. Thankfully, Mihalcea and Norton managed to link this reachability to a purely algebraic set $\mathcal{A}_d(u)$, which simplifies things considerably.

If we define $\mathcal{A}_d(u)$ as:
\begin{equation}
    \mathcal{A}_d(u) := \{ v \in D_\infty \mid \ell(uv) = \ell(u) + \ell(v) \text{ and } \varphi(v) \le d \},
\end{equation}
the bridge between building these chains and tracking the degree map $\varphi$ rests on a specific parity lemma. This lemma effectively dictates how degrees accumulate and stay bounded across chains:

\begin{lemma}[Parity of Chains, Lemma 2.5 in \cite{mihalcea2017}]
\label{lem:parity}
Take any $u, v \in D_\infty$ connected by a chain of degree $d'$. Then $d' = \varphi(u^{-1}v) + 2(r, s)$ for some integers $r, s \in \mathbb{Z}_{\ge 0}$. You'll notice this naturally forces $\varphi(u^{-1}v) \le d'$.
\end{lemma}

Armed with Lemma \ref{lem:parity}, we can put a hard bound on the degrees of any elements reachable from $u$. This builds up to Mihalcea and Norton's main theorem, which cleanly asserts that the curve neighborhood boils down to just tracking the maximal elements of $\mathcal{A}_d(u)$:
\begin{theorem}
\label{thm:main_math}
Let $u \in D_\infty$ and $d$ be a degree. Then the curve neighborhood is given by:
\begin{equation}
    \Gamma_d(u) = \{ uw \mid w \in \max \mathcal{A}_d(u) \}.
\end{equation}
\end{theorem}

Our core objective was to port Definition 1, Lemma \ref{lem:parity}, and Theorem \ref{thm:main_math} entirely into Lean. As you'll see in Section \ref{sec:core_proofs}, cementing Lemma \ref{lem:parity} natively pushes us to blend explicit modulo 2 arithmetic with the properties of Coxeter length—which makes it a prime candidate for machine proof.

\section{Survey of Lean Code}
\label{sec:formalization}

\subsection{Dihedral Groups and Coxeter Systems}
The infinite dihedral group \lstinline{DihedralGroup 0} is available in Mathlib~\cite{mathlib}. To leverage the API of Coxeter systems intuitively, we endow $D_\infty$ with a \lstinline{CoxeterSystem} instance by explicitly defining a matrix $M$ over \lstinline{Fin 2}.
\begin{lstlisting}[language=lean]
def M : CoxeterMatrix (Fin 2) := {
  M := !![1, 0; 0, 1]
  isSymm := by decide
  -- ... remaining fields proven by decide
}
\end{lstlisting}

The main engineering hurdle here isn't just stating the definitions; it's proving a clean multiplicative equivalence (\lstinline{MulEquiv}) between the geometric group $D_\infty$ and the abstract Coxeter group \lstinline{M.Group}. After defining two homomorphisms mapping the geometric reflections exactly onto $s_0$ and $s_1$, forcing Lean to realize they are identical inverses meant working through some intricate but tedious algebraic manipulation—specifically wrangling \lstinline{conj_eq_inv} to tame powers of alternating products. Once this was resolved, we locked in our Coxeter system:
\begin{lstlisting}[language=lean]
def cs : CoxeterSystem M D∞ := { mulEquiv := mulEquiv }
\end{lstlisting}
This setup grants us the right to cleanly write abstract Coxeter lengths as \lstinline{cs.length g} underneath the symbol $\ell(g)$.

\subsection{Explicit Lengths and Reduced Words}
One of the notable challenges when dealing with Coxeter groups formally is that \lstinline{cs.length} is defined non-constructively as a minimal list bound. Naturally, if we want to run computable operations to chase down curve neighborhoods, we really need a closed-form formula.

We bypassed this by constructing a computationally grounded \lstinline{explicit_length} function mirroring an explicit \lstinline{reducedWord} list (which itself leans on Mathlib's handy \lstinline{alternatingWord}). If you evaluate a positive reflection \lstinline{sr k}, for example, it generates an alternating sequence of length $2k - 1$:
\begin{lstlisting}[language=lean]
def explicit_length : D∞ → ℕ
| r k => 2 * k.natAbs
| sr k => if k > 0 then 2 * k.natAbs - 1 else 2 * k.natAbs + 1
\end{lstlisting}

This API is anchored entirely on proving \lstinline{length_eq}. Here, we verify that the mathematically pure but non-constructive abstract Coxeter length aligns exactly with this constructive function:
\begin{lstlisting}[language=lean]
theorem length_eq (g : D∞) : ℓ g = (reducedWord g).length
\end{lstlisting}
Getting this proof through relies on an inductive argument checking \lstinline{explicit_length g ≤ ℓ g}, and confirming \lstinline{reducedWord} safely maps back into $g$.

\subsection{Formalizing Degrees and the Moment Graph}
Degrees are tracked structurally with a basic \lstinline{Degree} object containing \lstinline{a b : ℕ}, complete with the necessary component-wise \lstinline{AddCommMonoid} ops and a \lstinline{PartialOrder}. When dealing with the edges forced by root reflections, we model them dynamically through predicates instead of building giant static lists:
\begin{lstlisting}[language=lean]
def IsEdge (u v : Vertex) (α : Root) : Prop := v = u * (s_α α)
\end{lstlisting}

Chains, meanwhile, manifest strictly as inductive sequences parameterized by endpoints and a running total degree:
\begin{lstlisting}[language=lean]
inductive HasIncreasingChain : Vertex → Vertex → Degree → Prop where
  | refl (u : Vertex) : HasIncreasingChain u u 0
  | step {u v w : Vertex} {d : Degree} {α : Root} :
      HasIncreasingChain u v d →
      IsEdge v w α →
      (ℓ v < ℓ w) →
      HasIncreasingChain u w (d + α.toDegree)
\end{lstlisting}

Through this structure, Bruhat traversal transforms into straightforward structural induction on chains, readily proving \lstinline{lt_iff_length_lt} that Bruhat order equates to length \cite[Lemma 2.3]{mihalcea2017}.

\section{Formalizing Core Lemmas}
\label{sec:core_proofs}

The mathematical justification heavily relies on reduced word combinatorics and parity preservation beneath the degree map. In the informal literature, the condition $\varphi(u^{-1}v) \le d'$ (Lemma \ref{lem:parity}) is established via implicit modulo 2 arithmetic over Coxeter relations. To formalize this seamlessly, we reduce the geometric parity constraints to explicit Presburger arithmetic in $\mathbb{N}$. Specifically, we formulate parity preservation as the existence of non-negative coefficients $r, s \in \mathbb{N}$:
\begin{lstlisting}[language=lean]
lemma degree_add_parity (g h : D∞) :
    ∃ (r s : ℕ), (φ g).a + (φ h).a = (φ (g * h)).a + 2 * r
               ∧ (φ g).b + (φ h).b = (φ (g * h)).b + 2 * s
\end{lstlisting}
Because evaluating the length and degree functions on $D_\infty$ introduces bounded additive perturbations, Lean's \lstinline{omega} tactic efficiently discharges the resulting linear arithmetic goals across all constructive cases of rotations and reflections. This algebraic foundation enables the structural induction required for Lemma \ref{lem:parity} (\lstinline{lemma_2_5_a} / \lstinline{lemma_2_5_b}).

\subsection{Finiteness and Main Theorem}
A critical hurdle in formalizing curve neighborhoods is establishing the finiteness of the algebraically defined set $\mathcal{A}_d(u)$. We overcome this by utilizing the explicit bounds $d_a + d_b + 1$ imposed by $\varphi$, effectively bounding the Coxeter length of any element $v \in \mathcal{A}_d(u)$. By systematically instantiating these bounds (\lstinline{h_len_bound}), we leverage Mathlib's \lstinline{exists_maximalFor} to  extract maximal elements within the Bruhat order:
\begin{lstlisting}[language=lean]
lemma exists_max_in_Ad (u : Vertex) (d : Degree) (z : Vertex) 
    (hz : z ∈ Ad u d) : ∃ w, IsMaximalIn w (Ad u d) ∧ z ≤ w
\end{lstlisting}

Connecting the abstract neighborhood $\Gamma_d(u)$ to this algebraic condition (Theorem \ref{thm:main_math}) requires reasoning about paths in the moment graph via Coxeter descents. We rely heavily on Mathlib's \lstinline{IsLeftDescent} API, coupled with structural manipulation of alternating words—via list operations like \lstinline{head?} and \lstinline{getLast?}—to construct descents dynamically. The proof fundamentally proceeds by length contradiction: if $\ell(v) < \ell(u) + \ell(u^{-1}v)$, we can extract a descent that yields a strictly shorter chain-preserving element, directly violating Bruhat maximality. This formally resolves the \lstinline{main_theorem}:
\begin{lstlisting}[language=lean]
theorem main_theorem (u : Vertex) (d : Degree) :
  CurveNeighborhood u d = {v | ∃ w, IsMaximalIn w (Ad u d) ∧ v = u * w}
\end{lstlisting}

\subsection{Computability}
\label{sec:computability}

In Lean 4, subsets defined purely mathematically—like \lstinline{Set α}—do not automatically provide an executable algorithm. To actually compute these curve neighborhoods, we had to carefully bridge the gap between classical verification and raw computation.

\textbf{Lengths and Executable Filtration.} First, abstract \lstinline{cs.length} was shadowed by a computationally tractable recursive function, \lstinline{cLength}, which we explicitly proved to be equivalent. Once paired with a \lstinline{Decidable} partial order on the degrees, generating bounded finite candidate sets suddenly became viable.

Knowing that elements in $\mathcal{A}_d(u)$ are strictly capped at lengths of up to $limit = d_a+d_b+1$, we could computationally enumerate the working elements straight into a \lstinline{Finset Vertex}:
\begin{lstlisting}[language=lean]
def enumerateD_list (n : ℕ) : List D∞ :=
  (List.range (n + 1)).flatMap fun k =>[cs.wordProd (alternatingWord 0 1 k), 
     cs.wordProd (alternatingWord 1 0 k)]

def enumerateD (n : ℕ) : Finset D∞ := (enumerateD_list n).toFinset
\end{lstlisting}

Filtering this \lstinline{enumerateD} set then yields an executable version of $\mathcal{A}_d(u)$:
\begin{lstlisting}[language=lean]
def Ad_finset (u : Vertex) (d : Degree) : Finset Vertex :=
  let limit := d.a + d.b + 1
  (enumerateD limit).filter (fun v => 
    cLength (u * v) = cLength u + cLength v ∧ φ v ≤ d)
\end{lstlisting}

From here, we extract the maximal elements and left-multiply by $u$, effectively producing a fully verified, strictly executable slice of logic that we call \lstinline{CurveNeighborhood_computable}.
\begin{lstlisting}[language=lean]
def CurveNeighborhood_computable (u : Vertex) (d : Degree) : 
    Finset Vertex :=
  let A := Ad_finset u d
  let maxA := A.filter (fun w => 
    ∀ w' ∈ A, ¬(cLength w < cLength w'))
  maxA.image (fun w => u * w)
\end{lstlisting}

\textbf{Native Evaluation.} Finally, leaning on Lean 4's \lstinline{#eval}, we can instantaneously reproduce Mihalcea and Norton's formulas \cite{mihalcea2017} directly within the dependent type system:
\begin{lstlisting}[language=lean]
-- Output: {r(2), r(-2)}
#eval CurveNeighborhood_computable 1 ⟨2, 2⟩ 

-- Output: {sr(-3)}
#eval CurveNeighborhood_computable s0 ⟨2, 3⟩
\end{lstlisting}
This setup tightly unites the formal mathematical proof with live, symbolic algebraic computation—all under one roof.

\section{Conclusion}
\label{sec:conclusion}

This formalization demonstrates that intricate combinatorial proofs in quantum Schubert calculus are fully accessible and computationally executable in Lean 4. Throughout this project, the robustness of Mathlib and the utility of instantiating $D_\infty$ as a \lstinline{CoxeterSystem} proved essential, effectively bridging the gap between explicit reduced word combinatorics and abstract algebraic hierarchies. 

A key insight from this work is the contrast between informal intuition and formal rigor. Several steps that are typically treated as immediate in the informal literature—such as the finiteness of $\mathcal{A}_d(u)$ or the preservation of modulo 2 parity—require significant technical effort to establish formally. In a machine-checked environment, these intuitive leaps must be replaced by explicit construction of bounds and exhaustive case analyzes. Despite the added complexity of this formal bookkeeping, by carefully separating low-level list manipulations from high-level geometric reasoning, we obtained a script that remains structurally isomorphic to the original arguments by Mihalcea and Norton.

Beyond theoretical verification, leveraging Lean's \lstinline{omega} tactic and \lstinline{Finset}-based representations enabled us to successfully unite formal proof with native computation. Ultimately, by establishing these moment graph results within the Lean ecosystem, we offer a reusable, mechanically verified framework for researchers to compute and verify exact algebraic neighborhoods with absolute certainty.

\section*{Acknowledgments}

This project was developed as part of the Lean Formalization course at the Academy of Mathematics and Systems Science, Chinese Academy of Sciences. We are deeply grateful to our instructor, Jujian Zhang, for running this course and for his constant encouragement throughout the development of this library. We would also like to thank the teaching assistant, Junqi Liu, for his dedicated Lean support. Many of the technical hurdles and proof strategies discussed in this paper were successfully resolved through the stimulating discussions held during the lectures.

\bibliographystyle{splncs04}
\bibliography{myrefs}

\end{document}